\documentclass[10pt,a4paper,twoside]{amsart}
\usepackage[utf8]{inputenc}
\usepackage[english]{babel}
\usepackage[T1]{fontenc}
\usepackage{microtype, fancyhdr, lmodern, xcolor}
\usepackage[a4paper, hmarginratio=1:1]{geometry}  
\usepackage{kerkis}
\usepackage{charter}

\usepackage{amsfonts, amsmath, amsthm, amssymb, mathrsfs, amscd, mathtools}
\usepackage[all]{xy}
\numberwithin{equation}{section}	
\usepackage{faktor}
\allowdisplaybreaks

\usepackage{url, enumitem}
\usepackage[noadjust]{cite}

\usepackage{multirow,bigdelim, caption}
\usepackage{booktabs}
\usepackage{graphicx}
\graphicspath{{Figures/}}

\usepackage{hyperref}

\theoremstyle{plain}
\newtheorem{theorem}{Theorem}[section]
\newtheorem{proposition}[theorem]{Proposition}
\newtheorem{corollary}[theorem]{Corollary}
\newtheorem{lemma}[theorem]{Lemma}

\theoremstyle{remark}
\newtheorem{remark}[theorem]{Remark}

\theoremstyle{definition}
\newtheorem{definition}[theorem]{Definition}

\newcommand\Zz{\mathbb{Z}}

\newcommand\Rr{\mathbb{R}}

\DeclareMathOperator{\Aut}{Aut}
\newcommand\BB[1]{B_{#1}}	
\newcommand\BC[2]{B_{#1}[#2]}	

\newcommand\inv{^{-1}} 
\newcommand\Mod[1]{\mathrm{Mod}({#1})} 	
\newcommand\PB[1]{P_{#1}}	

\newcommand\sig[1]{\sigma_{\hspace{-0.3ex}#1}^{\null}} 
\newcommand\sigg[2]{\sigma_{\hspace{-0.3ex}#1}^{#2}}	
\newcommand\Sp[2]{\mathrm{Sp}_{#1}(#2)}	
\newcommand\ii{i} 
\newcommand\nn{n} 
\newcommand\nno{n-1} 
\newcommand\mm{m} 

\begin{document}
\title{Congruence subgroups of braid groups and crystallographic quotients. Part I}

\author[Bellingeri]{Paolo Bellingeri}
\address{Normandie Univ, UNICAEN, CNRS, LMNO, 14000 Caen, France}
\email{paolo.bellingeri@unicaen.fr}

\author[Damiani]{Celeste Damiani}
\address{Fondazione Istituto Italiano di Tecnologia, Genova, Italy}
\email{celeste.damiani@iit.it}

\author[Ocampo]{Oscar Ocampo}
\address{Universidade Federal da Bahia, Departamento de Matem\'atica - IME, CEP: 40170 - 110, Salvador - BA, Brazil}
\email{oscaro@ufba.br}

\author[Stylianakis]{Charalampos Stylianakis}
\address{University of the Aegean, Department of mathematics, Karlovasi, 83200, Samos, Greece}
\email{stylianakisy2009@gmail.com}



\date{\today}

\begin{abstract}
This paper is the first of a two part series devoted to describing relations between congruence and crystallographic braid groups.  We recall and introduce some elements belonging to congruence braid groups and we establish some (iso)-morphisms between crystallographic braid groups
and corresponding quotients of congruence 
braid groups.
\end{abstract}

\maketitle

\section{Introduction}
This paper delves into the relationship between two families of groups, respectively subgroups and quotients of classical braid groups: \emph{congruence subgroups of braid groups} and \emph{crystallographic braid groups}, respectively introduced Arnol'd~\cite{Arnold:1968} and Tits~\cite{JTits:1966}.

While both families are instances of more general groups with rich theoretical backgrounds, they have also garnered significant attention in recent (and less recent) literature on braid groups and relatives, see for instance~\cite{Brendle-Margalit:2018, Stylianakis:2018, ABGH, Nakamura:2021, KordekMargalit:2022, bloomquist2023quotients} for congruence subgroups of braid groups and~\cite{ACampo:1979, Goncalves-Guaschi-Ocampo:2017, BeckMarin:2020, Goncalves-Guaschi-Ocampo-Pereiro, Bellingeri-Guaschi-Makri, Cerqueira-Ocampo} for crystallographic braid groups. Let us provide an overview of the two general families to which these groups belong.

In the context of groups of matrices, a \emph{congruence subgroup} of a matrix group with integer entries 
is a subgroup defined as the kernel of the mod $m$ reduction of a linear group.
The notion of congruence subgroups can be generalised for arithmetic subgroups of certain algebraic groups
for which we can define appropriate reduction maps.
A classical question about congruence subgroups
is the \emph{congruence subgroup problem}, first formulated in~\cite{Bass-Milnor-Serre:1967}: 
in this seminal paper Bass, Minor and Serre prove that for $n\geq 3$
the group $SL_\nn(\Zz)$ has the \emph{congruence subgroup property}, 
meaning that every finite-index subgroup of $SL_\nn(\Zz)$ 
contains a \emph{principal congruence subgroup}. 
The literature devoted to this problem in several settings is vast 
(we refer to~\cite{Raghunathan:2004} for a survey), 
linking the theory of arithmetic groups and geometric properties of related spaces.

In this spirit, we can define congruence subgroups of any group via a choice of representation
into~$\mathrm{GL}(n, \mathbb{Z})$. 
Let the braid group $\BB\nn$ be the mapping class group $\mathrm{Mod}(D_n)$ of 
the disc with $n$ marked points~$D_n$. We can define a symplectic representation 
and use it to define congruence subgroups of braid groups $\BC\nn\mm$. 
We will recall details  in Section~\ref{S:Congruence}, but let us give here an idea of the definition of these groups. We start with the integral Burau representation of $\BB\nn$, which 
is the representation $\rho \colon \BB\nn \rightarrow GL_\nn(\Zz)$ 
obtained by evaluating the (unreduced) Burau representation 
$\BB\nn \rightarrow GL_\nn(\Zz[t, t\inv])$ at~$t=-1$.
Describing the representation from a topological point of view, one can see that 
the integral Burau representation is symplectic, 
and can be regarded as a representation:
\begin{equation}
\label{E:rho}
\rho \colon \BB\nn \rightarrow 
\begin{cases} 
\Sp{\nn-1}{\Zz} \mbox{ for } \nn \mbox{ odd}, \\
(\Sp{\nn}{\Zz})_u \mbox{ for } \nn \mbox{ even}, \\
\end{cases}
\end{equation}
where $(\Sp {\nn} {\Zz} )_u$ is the subgroup of $\Sp \nn \Zz$ 
fixing a specific vector $u \in \Zz^\nn$, see~\cite[Proposition2.1]{Gambaudo-Ghys:2005} for a homological description of $(\Sp {\nn} {\Zz} )_u$ in this context.

The \emph{level $\mm$ congruence subgroup~$\BC\nn \mm$},
is the kernel of the mod $\mm$ reduction of the integral Burau representation
\begin{equation}
\label{E:rho_reduc}
\rho_\mm \colon \BB\nn \rightarrow
\begin{cases}
\Sp{\nn-1}{\Zz/\mm \Zz} \mbox{ for } \nn \mbox{ odd}, \\
(\Sp{\nn}{\Zz/\mm \Zz})_u \mbox{ for } \nn \mbox{ even}, \\
\end{cases}
\end{equation}
for~$\mm > 1$.

The second family of groups that we consider are \emph{crystallographic groups}, 
 appearing in the study of isometries of Euclidean spaces, see Section~\ref{S:Crystallographic}
 for precise definitions and useful characterisations. 
 In~\cite{Goncalves-Guaschi-Ocampo:2017}, Gon\c{c}alves, Guaschi and Ocampo prove
that certain quotients of the braid groups $\BB\nn$ are crystallographic, and use this result to study their torsion 
and other algebraic properties. 
The authors use this characterisation to prove
that the group $\faktor{\BB\nn}{[\PB\nn, \PB\nn]}$ is crystallographic, 
where $\PB\nn$ denotes the pure braid group on $\nn$ strands and $[\PB\nn, \PB\nn]$ its commutator subgroup. 
This quotient, that we will refer to as the \emph{crystallographic braid group}, was introduced by Tits in \cite{JTits:1966} as \emph{groupe de Coxeter
\'etendu}, see~\cite{Bellingeri-Guaschi-Makri} for a short survey.

Congruence subgroups and crystallographic structures share a point of contact. It follows from Arnol'd work~\cite{Arnold:1968}, that the pure braid group $\PB\nn$ can be characterised as the congruence subgroup $\BC\nn{2}$. With this equivalence and the results of~\cite{Goncalves-Guaschi-Ocampo:2017} in mind, it is natural to ask: how are congruence subgroups of braid groups and crystallographic groups related? This question was also recently raised in~\cite{Kumar:2024} for small Coxeter groups.
In this paper we propose to 
explore the interplay between congruence subgroups of braid groups and crystallographic groups, opening several questions that we will develop in a further work (\cite{BDMOS:2023}).
 
The paper is organised as follows. 
In Section~\ref{S:Congruence} we provide some basic definitions and properties that will be useful in this paper, such as the Burau representation, symplectic structures, the definition of congruence subgroups and the actions of half-twists on symplectic groups.
Section~\ref{S:Crystallographic} contains the main body of this work. In Subsection~\ref{Ss:3.1} we prove the following general result about crystallographic groups:
\newtheorem*{thm:mainthmcryst}{Theorem~\ref{mainthmcryst}}
\begin{thm:mainthmcryst}
Consider the short exact sequence
$
1 \longrightarrow K \longrightarrow G \stackrel{p}{\longrightarrow} Q \longrightarrow 1
$
where $K$ is a free abelian group of finite rank and $Q$ is a finite group such that the representation $\varphi\colon Q \to Aut(K)$, induced from the action by conjugacy, is not injective. 
Suppose that the group $p^{-1}(Ker(\varphi))$ is torsion free. 
Then $G$ is a crystallographic group with holonomy group $\faktor{Q}{Ker(\varphi)}$.
\end{thm:mainthmcryst}

This theorem will play an important role in this work, since the techniques used in \cite{Goncalves-Guaschi-Ocampo:2017} 
do not apply directly in this paper. This is because 
    the representation \[\Theta_m\colon \rho_m(\BB\nn) \to \mathrm{Aut}\left(\faktor{\BC\nn{m}}{[\BC\nn{m},\BC\nn{m}]}\right), \] induced from the action by conjugacy of $B_n$ on $\BC\nn{m}$, is injective if and only if $m =2$ (see Proposition~\ref{P:representation}), where $\rho_m$ is the homomorphism defined in \eqref{E:rho_reduc} .
We apply Theorem~\ref{mainthmcryst} to get the following result, that is proved in  
Subsection~\ref{Ss:3.2}, relating congruence subgroups and crystallographic groups.

\newtheorem*{thm:mainthmcongruenceb}{Theorem~\ref{mainthmcongruence}}
\begin{thm:mainthmcongruenceb}
    Let $n\geq 3$ be an odd integer and let $m\geq 3$ be a prime number. 
    If  the abelian group $\faktor{\BC\nn{m}}{[\BC\nn{m}, \BC\nn{m}]}$ is torsion free, then the group $\faktor{\BB\nn}{[\BC\nn{m}, \BC\nn{m}]}$ is crystallographic with dimension equal to rank$\left(\faktor{\BC\nn{m}}{[\BC\nn{m},\BC\nn{m}]}\right)$ and  holonomy group $\faktor{\rho_m(B_n)}{Z(\rho_m(B_n))}$. 
\end{thm:mainthmcongruenceb}

In Subsection~\ref{Ss:3.3} we show that there is an isomorphism between the crystallographic braid group $\faktor{\BB\nn}{[\PB\nn, \PB\nn]}$ and a quotient of congruence subgroups as described in the next result.

\newtheorem*{thm:second_main}{Theorem~\ref{T:second_main}}
\begin{thm:second_main}
Let $m$ be a positive integer and let $\nn \geq 3$. 
Consider the map \[ \overline{\xi}\colon \faktor{\BB\nn}{[\PB\nn, \PB\nn]}\to \faktor{\BC\nn{m}}{[\PB\nn, \PB\nn] \cap \BC\nn{m}}\] defined by $\overline{\xi}(\sigma_i)= \sigma_i^m$ for all $1\leq i\leq n-1$. 
If $m$ is odd, then $\overline{\xi}$ is an isomorphism. 
As a consequence, for $n\geq 3$ and $m$ odd,  $\faktor{\BC\nn{m}}{[\PB\nn, \PB\nn] \cap \BC\nn{m}}$ is a crystallographic group of dimension $n(n-1)/2$ and holonomy group $S_n$.
\end{thm:second_main}

Finally, in Appendix~\ref{S:Appendix} we give the proof of two technical lemmas that we use in this paper.

\subsection*{Acknowledgments}
The authors thank Karel Dekimpe for helpful conversations on crystallographic groups. We are also deeply thankful to Alan McLeay for invaluable contributions to the first versions of this work.
P.~B. was partially supported by the ANR project AlMaRe (ANR-19-CE40-0001). C.~D.  is a member of GNSAGA of INdAM, and was partially supported Leverhulme Trust research project grant “RPG-2018-029: Emergent Physics From Lattice Models of Higher Gauge Theory”.
O.~O. would like to thank Laboratoire de Math\'ematiques Nicolas Oresme (Universit\'e de Caen Normandie) for its hospitality from August 2023 to January 2024, where part of this project was developed and was partially supported by Capes/Programa Capes-Print/ Processo n\'umero 88887.835402/2023-00 and also by National Council for Scientific and Technological Development - CNPq through a \textit{Bolsa de Produtividade} 305422/2022-7.
The authors would like to thank the anonymous referee for careful reading and valuable comments. 

\section{Congruence subgroups}\label{S:Congruence}

Let $S$ be a connected, orientable surface, possibly with marked points and boundary components. The \emph{mapping class group} $\mathrm{Mod}(S)$ of $S$ is the group of homotopy classes of homeomorphisms of $S$ that preserve the orientation, fix the set of marked points setwise, and fix the boundary pointwise. 

\subsection{Braid groups and examples}\label{Subsec:braids}

\begin{figure}[htb]
	\centering
		\includegraphics[scale=.3]{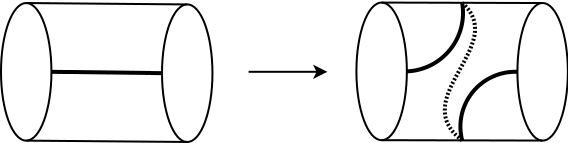}
	\caption{Twist map acts on an Annulus.}
	\label{twist}
\end{figure}

Let $S$ be surface as above. We introduce a particular element of $\mathrm{Mod}(S)$ that will be used throughout the paper. Let $A$ be an annulus. The homeomorphism depicted in Figure \ref{twist} is called a \emph{twist map}. 

Now, let $c \subset S$ be a simple closed curve. The regular neighborhood $\mathcal{N}(c)$ of $c$ is homeomorphic to an annulus $A$. Consider the homeomorphism $f_c$ that acts as a twist map on $\mathcal{N}(c)$ and as the identity on $S \setminus \mathcal{N}(c)$. The homotopy class of $f_c$ is called \emph{Dehn twist about $c$}, denoted by $T_c$ \cite[Section 3.1]{Farb-Margalit:book}.

Braid groups can be defined in several equivalent ways, long known to be equivalent, see for instance~\cite{Kassel:2008braid, BirmanBrendle:2005}). In this work it will be convenient to define them in terms of mapping class groups. Let $D_n$ be a disc with $n \in \mathbb{N}$ marked points in its interior. The \emph{braid group} $B_n$ is $\mathrm{Mod}(D_n)$.
For a geometric insight of twists in the context of braid groups, let $D_n$ lie on the $xy$-plane with its centre on the $x$-axis. Denote the punctures from left to right by $p_1, p_2, \ldots, p_n$: the arc connecting $p_i$ and $p_{i+1}$ is denoted by $a_i$ (see Figure \ref{cover}). Consider $a_i$ to be the diameter of a circle $c$ such that the points $p_i$ and $p_{i+1}$ lie on~$c$. Interchanging the points $p_i$ and $p_{i+1}$ by rotating them half way along $c$ in the clockwise direction gives a homeomorphism of $D_n$, and its homotopy class in  $\mathrm{Mod}(D_n)$ is called a \emph{half twist}, denoted by~$\sig\ii$. Note that all conjugates of $\sig\ii$ are called half twists. 
In terms of presented groups, half twists correspond to the Artin's generators from Artin's presentation for $\BB\nn$~\cite{Artin:1925}:
\[
\bigg\langle \sig1, \ldots\, , \sig\nno \ \bigg\vert \ 
\begin{matrix}
\sig{i} \sig j = \sig j \sig{i} 
&\text{for} &\vert  i-j\vert > 1\\
\sig{i} \sig j \sig{i} = \sig j \sig{i} \sig j 
&\text{for} &\vert  i-j\vert  = 1
\end{matrix}
\ \bigg\rangle.
\]
Let $c\in D_n$ be a curve surrounding the points $p_i, p_{i+1}$. This curve is homotopic to the circle described above. We note that if $\sig\ii$ is a half twist, then $\sigg\ii{2}$ is a Dehn twist about the curve $c$. This Dehn twist is generalised, for $1\leq i<j\leq n$, as follows
\[ A_{i,j} = (\sigma_{j-1} \sigma_{j-2} \ldots \sigma_{i+1}) \sigg{i}{2} (\sigma_{j-1} \sigma_{j-2} \ldots \sigma_{i+1})^{-1}. \]

\begin{figure}[htb]
	\centering
		\includegraphics[scale=.45]{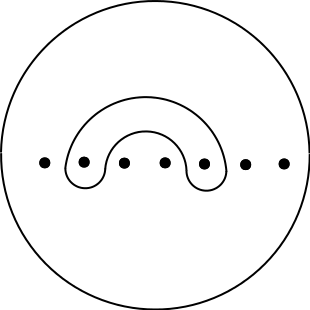}
	\caption{The Dehn twist along the curve that surrounds the punctures $p_2,p_5$ is $A_{2,5}$.}
	\label{a_ij}
\end{figure}

We recall that a generating set of $P_{n}$ is given by $\{A_{i,j}\}_{1\leq i<j\leq n}$. 
Geometrically the element $A_{i,j}$ can be represented as a Dehn twist about a curve surrounding punctures $p_i$ and $p_j$. For instance in Figure~\ref{a_ij} we describe $A_{2,5}$ as the Dehn twist about the curve that surrounds punctures $p_2$ and $p_5$. 

We shall be interested in the action by conjugation of $B_n$ on $P_n$. 
Recall from  \cite[Proposition~3.7, Chapter 3]{Murasugi} that for all $1\leq k\leq n-1$ and for all $1\leq i<j\leq n$,
\begin{equation*}\label{eq:conjugAij}
\sigma_kA_{i,j}\sigma_k^{-1}=
\begin{cases}
A_{i,j} & \text{if $k\neq i-1, i, j-1, j$}\\
A_{i,j+1} & \text{if $j=k$}\\
A_{i,j}^{-1}A_{i,j-1}A_{i,j} & \text{if $j=k+1$ and $i<k$}\\
A_{i,j} & \text{if $j=k+1$ and $i=k$}\\
A_{i+1,j} & \text{if $i=k<j-1$}\\
A_{i,j}^{-1}A_{i-1,j}A_{i,j} & \text{if $i=k+1$}.
\end{cases}
\end{equation*}
This action induces an action of $\faktor{B_n}{[P_n, P_n]}$ on $\faktor{P_n}{[P_n, P_n]}$, see \cite[Proposition~12]{Goncalves-Guaschi-Ocampo:2017}: Let $\alpha\in \faktor{B_n}{[P_n, P_n]}$, and let $\pi$ be the permutation induced by $\alpha^{-1}$, then $\alpha A_{i,j}\alpha^{-1}=A_{\pi(i), \pi(j)}$ in $\faktor{P_n}{[P_n, P_n]}$.

Another important element of $B_n$ that will play crucial role in the paper is the Dehn twist (or a full twist) along a curve surrounding all marked points of $D_n$. We denote this element by $\Delta_n^2$. In fact $\Delta_n^2$ generates the center of $B_n$ \cite{Chow}; in terms of half twists, we have
\[ \Delta_n^2 = (\sig1 \sig2 \ldots \sig\nno)^n. \]

\subsection{Burau representation and symplectic structures}
\label{Ss:2.2}

Braid groups naturally act on the homology of topological spaces obtained from the punctured disk. A construction arising in such a way is the Burau representation~\cite{Burau:1935}.
One of the most famous representations of the braid group, originally introduced in terms of matrices assigned to the generators in the Artin's
presentation of $\BB\nn$, the Burau representation is fundamental in low-dimensional topology. While this representation has been extensively studied, it still retains some mystery:
a long standing candidate for proving the linearity of the braid group (later established independently in~\cite{Bigelow:2001} and~\cite{Krammer:2002}), the question of its faithfulness  has remained open for quite some time. The Burau representation, faithful for $\nn\leq 3$~\cite{MagnusPeluso:1969}, eventually proved to be unfaithful for $\nn \geq 5$ (Moody~\cite{Moody:1991} proved unfaithfulness for $\nn \geq 9$, Long and Paton~\cite{LongPaton:1993} for $\nn \geq 6$, and Bigelow~\cite{Bigelow:1999} for $\nn = 5$). However, the case $\nn=4$ remains open, with advances towards closing the problem being published recently~\cite{Traczyk:2018, Birman:2021, datta:2022}.

In this work we are going to take the viewpoint of the Burau representation as a homological representation. Let $\pi=\pi_1(D_n,q)$ denote the fundamental group of $D_n$ where $q \in \partial D_n$. The function  $\pi \to \Zz\cong \langle t   \rangle$ defines a covering space $\Tilde{D}_n \to D_n$. Let $Q$ be a set of all lifts of $q$. The action of $t$ on $\Tilde{D}_n$ induces a $\Zz[t]$-module $\mathrm{H}_1(\Tilde{D}_n,Q;\Zz[t])$ of dimension $n$. Every mapping class in $\mathrm{Mod}(D_n)$ lifts to a unique mapping class in $\Tilde{D}_n$. Hence, the (reducible) Burau representation is given by a map
\[
 \mathrm{Mod}(D_n) \to \mathrm{Aut}(\mathrm{H}_1(\Tilde{D}_n,Q;\Zz[t])).
\]
This representation splits into a direct sum of an $(n-1)$ and a $1$ dimensional representations.

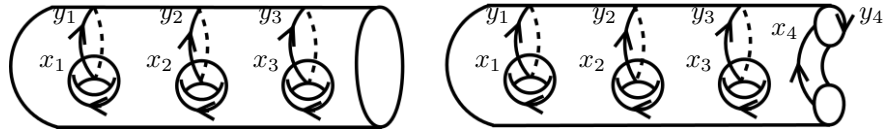
\begin{figure}[htb]
	\centering
		\includegraphics[scale=.45]{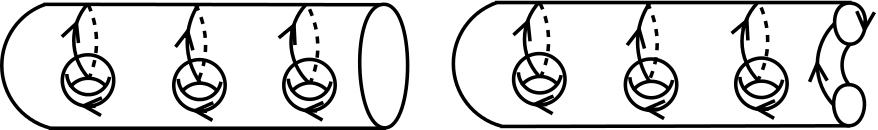}
	\caption{Generators for $\mathrm{H}_1(\Sigma';\Zz)$ on the left, and $\mathrm{H}_1(\Sigma,Q;\Zz)$ on the right.}
	\begin{picture}(22,12)
		\put(-130,78){$y_1$}
		\put(-135,60){$x_1$}
		\put(-90,78){$y_2$}
		\put(-95,60){$x_2$}
		\put(-53,78){$y_3$}
		\put(-55,60){$x_3$}
		\put(32,78){$y_1$}
		\put(28,60){$x_1$}
		\put(72,78){$y_2$}
		\put(67,60){$x_2$}
		\put(109,78){$y_3$}
		\put(107,60){$x_3$}
		\put(172,78){$y_4$}
		\put(139,72){$x_4$}
	\end{picture}
	\label{sympb1}
\end{figure}

Fixing $t=-1$, the covering space becomes a two-fold branch cover $\Sigma \to D_n$, where $\Sigma$ is homeomorphic to a surface of genus $g=(n-1)/2$ and one boundary component if $n$ is odd, and $g=n/2 - 1$ and two boundary components if $n$ is even  \cite{PerronVannier:1996}. As mentioned above, every mapping class in $\mathrm{Mod}(D_n)$ lifts to a unique mapping class in $\mathrm{Mod}(\Sigma)$ leading to an injection $\mathrm{Mod}(D_n) \to \mathrm{Mod}(\Sigma)$. Let $q \in \partial D_n$ be a point and $Q$ be a set of all lifts of $q$. The reducible Burau representation at $t=-1$ \cite[Section 2]{Brendle-Margalit:2018} (see also \cite{bloomquist2023quotients}) is 
\[
\mathrm{Mod}(D_n) \to \mathrm{Mod}(\Sigma) \to \mathrm{Aut}(\mathrm{H}_1(\Sigma, Q;\Zz)).
\] 
For $n$ odd, the module $\mathrm{H}_1(\Sigma, Q;\Zz)$ splits as $\mathrm{H}_1(\Sigma;\Zz) \times \Zz$ and 
the induced action of $\mathrm{Mod}(D_n)$ preserves a symplectic form on $\mathrm{H}_1(\Sigma;\Zz)$. Hence, the image of the latter representation is conjugate to $\Sp\nno\Zz$~\cite[Proposition~2.1]{Gambaudo-Ghys:2005}.
When $n$ is even, the module $\mathrm{H}_1(\Sigma,Q;\Zz)$ carries a symplectic structure. More precisely, if $g$ is the genus of $\Sigma$, then let $\Sigma'$ be a surface obtained by gluing a pair of pants in the boundary of $\Sigma$. Then $\Sigma'$ is a surface genus $g+1$ with one boundary component. We consider $\mathrm{H}_1(\Sigma,Q;\Zz)$ as a submodule of $\mathrm{H}_1(\Sigma';\Zz)$. In Figure \ref{sympb1} we give a basis for each of the latter modules.

The representation obtained by the construction above is

\begin{equation}
\label{E:rho_reduc_z}
\rho \colon \BB\nn \rightarrow 
\begin{cases} 
\Sp{\nn-1}{\Zz} \mbox{ for } \nn \mbox{ odd}, \\
(\Sp{\nn}{\Zz})_u \mbox{ for } \nn \mbox{ even}, \\
\end{cases}
\end{equation}
where, without loss of generality, we can choose $u = y_{2g+1}$.
For the detailed construction, 
see~\cite[Section~2.1]{Brendle-Margalit:2018}. 

An analogue of the principal congruence subgroups 
for the braid groups $\BB\nn$ can be defined 
starting from integral Burau representation. 
The \emph{level $\mm$ congruence subgroup~$\BC\nn \mm$}
is the kernel of the mod $\mm$ reduction of the integral Burau representation
\begin{equation} 
\rho_\mm \colon \BB\nn \rightarrow 
\begin{cases} 
\Sp{\nn-1}{\Zz/\mm \Zz} \mbox{ for } \nn \mbox{ odd}, \\
(\Sp{\nn}{\Zz/\mm \Zz})_u \mbox{ for } \nn \mbox{ even}, \\
\end{cases}
\end{equation}
for~$\mm > 1$.

In~\cite{Arnold:1968}, Arnol'd proved that the pure braid group~$\PB\nn$ is isomorphic to the level $2$ congruence subgroup $\BC\nn 2$ of the braid group $\BB\nn$, see also~\cite[Section~2]{Brendle-Margalit:2018} for a sketch of the original argument. 
In \cite{Brendle-Margalit:2018}, Brendle and Margalit go on to prove that $\BC\nn 4$ is isomorphic to the subgroup~$\PB\nn^2$, where $\PB\nn^2$
is the subgroup of $\PB\nn$ generated
by the squares of all elements.

A well-known family of elements in $\BC\nn \mm$
are \emph{Braid Torelli elements}. Consider the symplectic representation~\eqref{E:rho_reduc_z}. The kernel of this representation is denoted by $\mathcal{BI}_n$ and it is called \emph{braid Torelli} group. Since the representation \eqref{E:rho_reduc} is a $\bmod(m)$ reduction of $\rho$, then every element of $\mathcal{BI}_n$ is actually an element of $\BC\nn{m}$. In particular, $\mathcal{BI}_n$ is generated by squares of Dehn twists about curves surrounding odd number of marked points in $D_n$ \cite{BMP:2015}. In terms of half-twists, these elements are of the form
 \[ (\sig1 \ldots \sigma_k)^{2k+2} ,\]
where $k < n$ is even. This family of elements can be extended. If for example we denote by $c$ a curve surrounding odd number of marked points, then $T_{c}^2 \in \mathcal{BI}_n$). Other families of elements
in $\BC\nn \mm$, like as mod $p$ involutions and center maps are described in \cite[Section 4]{Stylianakis:2018}.

\subsection{Actions of half-twists on symplectic groups}

Recall that $B_n \cong \mathrm{Mod}(D_n)$ and $\Sigma \to D_n$ be a two-fold branched cover. The image of the monomorphism $\mathrm{Mod}(D_n) \to \mathrm{Mod}(\Sigma)$ is called the \emph{hyperelliptic mapping class group} denoted by $\mathrm{SMod}(\Sigma)$. Below we explain how to lift elements of $\mathrm{Mod}(D_n)$ into $\mathrm{SMod}(\Sigma)$. Then we use these lifts to explain their action on $\mathrm{H}_1(\Sigma,Q;\Zz)$.

\begin{figure}[htb]
	\centering
		\includegraphics[scale=.45]{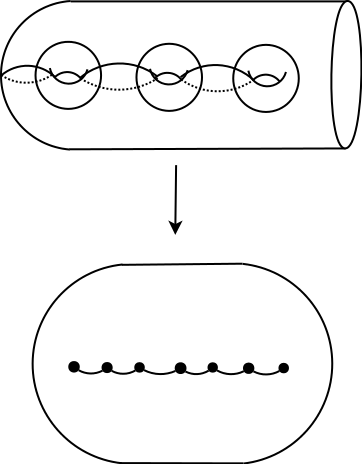}
	\caption{An example of a 2-fold cover of a marked disc. The simple closed curve $c_i$ in the genus 3 surface becomes the arc $a_i$ in the disc.}
	\begin{picture}(22,12)
		\put(-33,85){$a_1$}
            \put(-18,85){$a_2$}
            \put(0,85){$a_3$}
            \put(18,85){$a_4$}
            \put(33,85){$a_5$}
            \put(50,85){$a_6$}
		\put(-64,215){$c_1$}
            \put(-20,213){$c_3$}
            \put(25,210){$c_5$}
            \put(-44,205){$c_2$}
            \put(0,203){$c_4$}
            \put(45,203){$c_6$}
	\end{picture}
	\label{cover}
\end{figure}
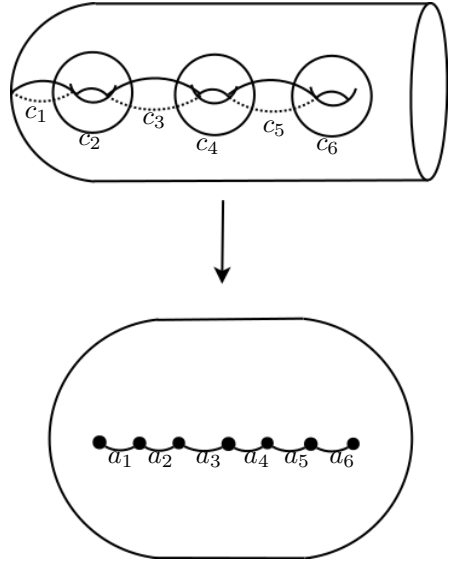

Let $\Sigma$ be a genus $g$ surface as in Figure \ref{cover}. The surface $\Sigma$ is the 2-fold cover of the disc $D_n$. Each simple closed curve $c_i$ is a lift of the arc $a_i$. Recall that $\sig\ii$ is a half twist along $a_i$. Then $\sig\ii$ lifts to the Dehn twist $T_{c_i}$. This association describes the homomorphism $\BB\nn \to \mathrm{SMod}(\Sigma)$ by $\sig\ii \mapsto T_{c_i}$.

Suppose that $\Sigma$ is a genus $g \geq 1$ surface with one boundary component (similarly for two boundary components). Let $T_c$ be a Dehn twist about a simple closed curve $c$ and let $[c]$ be its homology class in $\mathrm{H}_1(\Sigma;\Zz)$. Denote by $t_{[c]}$ a transvection induced by $T_c$. The action of the transvection $t_{[c]}$ on a homology class $u$ is defined by $t_{[c]}(u) = u + i(u,[c])[c]$ where $i(,)$ is a symplectic form. Therefore, the homomorphism $\rho_m \colon B_n \to \Sp{\nn-1}{\Zz/\mm \Zz}$ is defined by $\sigma_i \mapsto t_{[c_i]}$ (similarly for two boundary components). The next two lemmas describe the image of particular elements of $\BB\nn$ to the symplectic group over $\Zz/\mm \Zz$.

\begin{lemma}
For $m\geq 2 $ we have that $\rho_m(\sigma^m_i)=1$.
\label{m_powers}
\end{lemma}

\proof
Since $\sig\ii$ is mapped to the transvection $t_{[c_i]}$, we only need to compute the matrix form of $t_{[c_i]}$. It is easy to calculate the action of $t_{[c_i]}$ on the basis of Figure \ref{sympb1}. The result is conjugate to the following matrix:
\[
\begin{pmatrix}
1&1\\	
0&1
\end{pmatrix} \oplus I,
\]
where $I$ is the identity matrix of dimension $n-2$. The result follows by calculating the $m$-th power of the latter matrix over $\Zz/ \mm \Zz$.
\endproof

Lemma \ref{m_powers} leads to the question if $\BC\nn \mm$ coincides
with the group normally generated by $\sigma_i^m$. This is generally not the case (see 
\cite{BDMOS:2023} for further details): 
in fact $\BC\nn \mm$ is of finite index in $B_n$ while the group normally generated by $\sigma_i^m$
is not (except pairs $(n,m)\in\{ (3,3), (3,4), (3,5), (4,3), (5,3) \}$, see
\cite{Coxeter:1957}).

Recall that  $\Delta_n^2$ denotes the element $(\sigma_1 \sigma_2 \ldots \sigma_{n-1})^{n}$ in $\BB\nn$, generating the center of $\BB\nn$.

\begin{remark}
The full twist $\Delta_n^2$ has this notation since it is the square of the Garside element $\Delta_n$, which is another crucial element in braid theory.    
\end{remark}

\begin{lemma}
If $n$ is odd, then $\rho_m(\Delta_n^2)$ has order 2. If $n$ is even, then $\rho_m(\Delta_n^2)$ has order $m$ if $gcd(2, m) = 1$ or it has order $m / 2$ if $gcd(2,m)=2$.
\label{center_image}
\end{lemma}

\proof
Suppose that $n$ is odd. The lift of $(\sigma_1 \sigma_2 \ldots \sigma_{n-1})^{n}$ to $\Sigma$ is the product of Dehn twists $(T_{c_1} T_{c_2} \ldots T_{c_{n-1}})^{n}$. Consider the basis $\{x_i, y_i \}$ depicted in Figure \ref{sympb1}. Then the action of the product $(t_{[c_1]} t_{[c_2]} \ldots t_{[c_{n-1}]})^{n}$ reverses the orientation of $x_i, y_i$ \cite{Stylianakis:2018}. 
Thus, it has order 2.

Suppose that $n$ is even. The lift of $(\sigma_1 \sigma_2 \ldots \sigma_{n-1})^n$ to $\Sigma$ is the product $(T_{c_1} T_{c_2} \ldots T_{c_{n-1}})^{n}$. By the chain relation, the latter product is $T_{q_1} T_{q_2}$ where the curves $q_1,q_2$ are parallel to the boundary components of $\Sigma$ \cite[Proposition 4.12]{Farb-Margalit:book}. Since $[q_1]=[q_2] = y_{n-1}$, we have that $T_{q_1} T_{q_2}$ is mapped into the square transvection $t_{y_{n-1}}^2$. The transvection $t_{y_{n-1}}$ fixes all basis elements $\{ x_i, y_i\}$ except $x_{n-1}$. Hence, 
\[ t_{y_{n-1}}^2(x_{n-1}) = x_{n-1} + 2 y_{n-1}. \]
\endproof

\section{Crystallographic structures and congruence subgroups of the braid groups}\label{S:Crystallographic}

We recall the definition of a crystallographic group.

\begin{definition}
A group $G$ is said to be a \emph{crystallographic group} 
if it is a discrete and uniform subgroup 
of $\Rr^N \rtimes \mathrm{O}(N, \Rr) \subseteq \mathrm{Aff}(\Rr^N)$.
\end{definition}

In~\cite{Goncalves-Guaschi-Ocampo:2017} there is a characterisation of crystallographic groups that is convenient in our context, see also \cite[Section~2.1]{Dekimpe}.  \

\begin{lemma}[{\cite[Lemma~8]{Goncalves-Guaschi-Ocampo:2017}}] 
\label{L:crystallo_char} 
A group $G$ is crystallographic if and only if there is 
and integer $N$ and a short exact sequence
\[
\begin{CD}
1 @>>> \Zz^N @>>> G  @>\zeta>> \Phi  @>>> 1
\end{CD}
\]
such that 
\begin{enumerate}
\item $\Phi$ is finite;
\item the integral representation  
$\Theta \colon \Phi \rightarrow \mathrm{Aut}(\Zz^N)$, 
induced by conjugation on $\Zz^N$ and defined by 
$\Theta(\phi)(x)= \pi x \pi\inv$, where $x \in \Zz^N$, $\phi \in \Theta$ 
and $\pi \in G$ is such that $\zeta(\pi)=\phi$ is faithful.
\end{enumerate}
\end{lemma}

\subsection{A general result on crystallographic groups}\label{Ss:3.1}

We prove in this subsection two results that are general, and that will be applied to the study of crystallographic structures on quotients of the braid group by commutator subgroups of congruence subgroups.

\begin{theorem}
    Let $\phi\colon G\to F$ be a surjective homomorphism with $F$ a finite group. 
    Let $K$ denote the kernel of $\phi$. Suppose that there is a non-trivial element of the center of $G$ that does not belong to $K$.
    Then the representation $\eta\colon F\to \mathrm{Aut}\left(\faktor{K}{[K,K]}\right)$, induced from the action by conjugacy of $\faktor{G}{[K,K]}$ on $\faktor{K}{[K,K]}$,  is not injective.
    \label{key_theorem_non_injective}
\end{theorem}

\proof 
Since $[K,K]$ is characteristic in $K$ and $K$ is normal in $G$, then $[K,K]$ is normal in $G$. 
Hence, we may consider the action by conjugacy of $\faktor{G}{[K,K]}$ on $\faktor{K}{[K,K]}$. 
This induces a representation $\eta\colon F\to \mathrm{Aut}\left(\faktor{K}{[K,K]}\right)$.
Let $z\in Z(G)$ be a non-trivial element in the center of $G$ such that $z\notin K$. 
We note that $\overline{z}$ does not belong to $\faktor{K}{[K,K]}$. Furthermore, since  $z\in Z(G)$, then 
\begin{equation}\label{eqncentre}
\overline{z}\overline{k}\overline{z}^{-1}=\overline{k}, \textrm{ for every element $\overline{k}\in \faktor{K}{[K,K]}$.}
\end{equation}
Let $\overline{\phi}(\overline{z})=t$, where $\overline{\phi}\colon \faktor{G}{[K,K]}\to F$. Notice that $t$ is a non-trivial element in $F$. 
So, we conclude that $\eta$ is not injective since $\eta(t)$ is the identity homomorphism (see \eqref{eqncentre}).
\endproof

In the following result, we consider the case where

the holonomy representation defined in Lemma~\ref{L:crystallo_char} is not injective and give conditions for the middle group to be a crystallographic group.

\begin{theorem}    
Consider the short exact sequence
$
1 \longrightarrow K \longrightarrow G \stackrel{p}{\longrightarrow} Q \longrightarrow 1
$
where $K$ is a free abelian group of finite rank and $Q$ is a finite group such that the representation $\varphi\colon Q \to Aut(K)$, induced from the action by conjugacy, is not injective. 
Suppose that the group $p^{-1}(Ker(\varphi))$ is torsion free. 
Then $G$ is a crystallographic group with holonomy group $\faktor{Q}{Ker(\varphi)}$.
\label{mainthmcryst}
\end{theorem}

\proof

First we note that $p^{-1}(Ker(\varphi))$ is a Bieberbach group, since it is finitely generated, torsion free and virtually abelian, see \cite[Theorem~3.1.3~(4)]{Dekimpe}.

Now, we prove that $p^{-1}(Ker(\varphi))$ is free abelian.
Since $p^{-1}(Ker(\varphi))$ is a Bieberbach group, then it fits in a short exact sequence
$
1 \to A \to p^{-1}(Ker (\varphi) ) \to F \to 1
$
where $F$ is a finite group, and $A$ is a free abelian group containing $K$ as
a normal subgroup of finite index. Suppose now that $F$ is not the trivial group.
Let $x \in p^{-1}(Ker (\varphi))$ be an element which is mapped onto a non
trivial element in $F$. We know that the induced map $F \to Aut(A)$ is injective,
so conjugation by $x$ induces a nontrivial automorphism of $A$. But since $K$ is of finite index in the free abelian group $A$, this implies that
conjugation by $x$
also induces a nontrivial automorphism of $K$. But this is not possible
since $x \in p^{-1}(Ker (\varphi) )$.

 Hence, $p^{-1}(Ker(\varphi))$ is free abelian, and we obtain the sequence
$
1 \longrightarrow p^{-1}(Ker(\varphi)) \longrightarrow G \stackrel{\overline{p}}{\longrightarrow} \faktor{Q}{Ker(\varphi)} \longrightarrow 1
$
such that the middle group is a crystallographic group. 
\endproof

\subsection{Crystallographic structures and congruence subgroups of braid groups}
\label{Ss:3.2}
In this subsection we study a quotient of $\BB\nn$, namely, $\faktor{\BB\nn}{[\BC\nn{m}, \BC\nn{m}]}$. Since $\faktor{\BB\nn}{[\BC\nn{2}, \BC\nn{2}]}$ is crystallographic \cite[Proposition~1]{Goncalves-Guaschi-Ocampo:2017}, being isomorphic to the crystallographic braid group $\faktor{\BB\nn}{[\PB\nn, \PB\nn]}$,  it is reasonable to ask whether $\faktor{\BB\nn}{[\BC\nn{m}, \BC\nn{m}]}$ is crystallographic for any positive integer $m$. 
Here we give conditions for this statement to holds.

The following short exact sequence
\[
\begin{CD}
1 @>>> \BC\nn{m} @>>> \BB\nn @>{\rho_m}>> \rho_m(\BB\nn) @>>> 1.
\end{CD}
\]
induces a short exact sequence on the quotients
\begin{equation}\label{E:rhom}
\begin{CD}
1 @>>> \faktor{\BC\nn{m}}{[\BC\nn{m},\BC\nn{m}]} @>>> \faktor{\BB\nn}{[\BC\nn{m},\BC\nn{m}]} @>{\overline{\pi}}>>  \rho_m(\BB\nn) @>>> 1.
\end{CD}
\end{equation}

The action by conjugacy of $\faktor{\BB\nn}{[\BC\nn{m},\BC\nn{m}]}$ on $\faktor{\BC\nn{m}}{[\BC\nn{m},\BC\nn{m}]}$ induces a homomorphism 
\begin{equation}\label{E:thetam}
\Theta_m\colon \rho_m(\BB\nn) \to \Aut\left(\faktor{\BC\nn{m}}{[\BC\nn{m},\BC\nn{m}]}\right). 
\end{equation}

As a consequence of Theorem~\ref{key_theorem_non_injective} we have the following result.

\begin{proposition}
\label{P:representation}
    The representation $\Theta_m\colon \rho_m(\BB\nn) \to \mathrm{Aut}\left(\faktor{\BC\nn{m}}{[\BC\nn{m},\BC\nn{m}]}\right)$, induced from the action by conjugacy of $B_n$ on $\BC\nn{m}$, is injective if and only if $m =2$.
\end{proposition}

\proof
For $m=2$, the abelian group $\faktor{\BC\nn{2}}{[\BC\nn{2},\BC\nn{2}]}$ has finite rank and it is torsion free. Furthermore, $\Theta_2$ in injective, see \cite[Proof of Proposition~1]{Goncalves-Guaschi-Ocampo:2017}.

Let $m\geq 3$. 
Recall that the element $\Delta_n^2 = (\sig1 \sig2 \ldots \sig\nno)^{n}$ represents the full twist on $\Mod{D_n} \cong \BB\nn$, which generates the center of $B_n$. 
From Lemma~\ref{center_image}, for any $n$, the element $\rho_m(\Delta^2_n)$ is non trivial and of finite order.
Thus, $\Delta^2_n \notin \BC\nn{m}$. Therefore, the induced element in $\faktor{\BB\nn}{[\BC\nn{m},\BC\nn{m}]}$ does not belong to $\faktor{\BC\nn{m}}{[\BC\nn{m},\BC\nn{m}]}$. 
From Theorem \ref{key_theorem_non_injective}, the homomorphism $\Theta_m\colon \rho_m(\BB\nn) \to \mathrm{Aut}\left(\faktor{\BC\nn{m}}{[\BC\nn{m},\BC\nn{m}]}\right)$ is not injective.
\endproof

Since the representation $\Theta_m$ is not injective for $m\geq 3$ we cannot apply Lemma~\ref{L:crystallo_char} in this case. 
However, we may give general conditions such that the group $\faktor{\BB\nn}{[\BC\nn{m},\BC\nn{m}]}$ is crystallographic. 
We have the following result about crystallographic structures and quotients of braid groups by commutators of congruence subgroups.

\begin{theorem}
    Let $n\geq 3$ be an odd integer and let $m\geq 3$ be a prime number. 
    If  the abelian group $\faktor{\BC\nn{m}}{[\BC\nn{m}, \BC\nn{m}]}$ is torsion free, then the group $\faktor{\BB\nn}{[\BC\nn{m}, \BC\nn{m}]}$ is crystallographic with dimension equal to rank$\left(\faktor{\BC\nn{m}}{[\BC\nn{m},\BC\nn{m}]}\right)$ and  holonomy group $\faktor{\rho_m(B_n)}{Z(\rho_m(B_n))}$. 
 \label{mainthmcongruence}
\end{theorem}

\proof
From Theorem~\ref{mainthmcryst} if $\overline{\rho_m}^{-1}(Ker(\Theta_m))$ is torsion free, where $\rho_m$ and $\Theta_m$ are the homomorphisms defined in \eqref{E:rhom} and \eqref{E:thetam}, respectively, then the group $\faktor{\BB\nn}{[\BC\nn{m}, \BC\nn{m}]}$ is crystallographic.

We note that the $Ker(\Theta_m)$ is isomorphic to $Z(\rho_m(B_n))$ the center of $\rho_m(B_n)$, since the full twist $\Delta_n^2$ generates the center of $B_n$ and $\rho_m(\Delta_n^2)$ belongs to the normal subgroup  $Ker(\Theta_m)$ of the symplectic group $\rho_m(B_n)$.
Recall that, under the assumptions of the statement, $Z(\rho_m(B_n))$ is isomorphic to $\Zz/2\Zz$. 
We consider now the following short exact sequence
\[
\begin{CD}
1 @>>> \faktor{\BC\nn{m}}{[\BC\nn{m},\BC\nn{m}]} @>>> \overline{\rho_m}^{-1}(Ker(\Theta_m)) @>{\overline{\rho_m}}>> Z(\rho_m(B_n)) @>>> 1
\end{CD}
\]
such that the kernel is torsion free (by hypothesis), the class of the element $\Delta_n^2\in B_n$ is a non-trivial element of $\overline{\rho_m}^{-1}(Ker(\Theta_m))$ and $1\neq \overline{\Delta_n^4}\in \faktor{\BC\nn{m}}{[\BC\nn{m},\BC\nn{m}]}$.

Applying a standard method to give presentations for group extensions~\cite[Chapter~10]{Johnson} and using the fact that the full twist generates the center of $B_n$, we conclude that the middle group $\overline{\rho_m}^{-1}(Ker(\Theta_m))$ is free abelian, and its rank corresponds to the rank of the free abelian group $\faktor{\BC\nn{m}}{[\BC\nn{m},\BC\nn{m}]}$.
\endproof

\begin{remark}

As far as we know, it is still an open problem whether $\faktor{\BC\nn{m}}{[\BC\nn{m},\BC\nn{m}]}$ is torsion free for any $n$ and $m$ except few cases. 
It is well known that the group $\faktor{\BC\nn{2}}{[\BC\nn{2},\BC\nn{2}]}$ is free abelian of rank $\binom{n}{2}$. 
Also, the groups $\faktor{\BC3{3}}{[\BC3{3},\BC3{3}]}$ and $\faktor{\BC3{4}}{[\BC3{4},\BC3{4}]}$ are torsion free of rank $4$ and 6, respectively, see  \cite{BDMOS:2023}.
\end{remark}

\subsection{Symmetric quotients of congruence subgroups of braid groups}
\label{Ss:3.3}

From the definition of congruence subgroups we get an inclusion $\iota\colon B_n[m]\to B_n$ that induces a homomorphism $\overline{\iota}\colon \faktor{B_n[m]}{[P_n, P_n]\cap B_n[m]}\to \faktor{B_n}{[P_n, P_n]}$. 
In general, $\overline{\iota}$ is not an isomorphism. In the following result we study it in more detail.

\begin{theorem}
\label{T:first_main}
Let $m$ be an odd positive integer and let $\nn \geq 3$.
The homomorphism induced from the inclusion $\iota\colon B_n[m]\to B_n$
\[
\overline{\iota}\colon \faktor{B_n[m]}{[P_n, P_n]\cap B_n[m]}\to \faktor{B_n}{[P_n, P_n]}
\]
is injective. 
Furthermore, the group $\overline{\iota}\left(\faktor{\BC\nn{m}}{[\PB\nn, \PB\nn] \cap \BC\nn{m}}\right)$ is a normal proper subgroup of  $\faktor{\BB\nn}{[\PB\nn, \PB\nn]}$ such that the quotient is isomorphic to $(\Zz/\mm \Zz)^{n(n-1)/2}$.
\end{theorem}

\begin{remark}
    For $n=2$ the quotient groups of Theorem~\ref{T:first_main} are isomorphic.
\end{remark}
Before delving into the proof, we state two technical lemmas that are going to be needed. However, we defer the proofs of said lemmas to Appendix~\ref{S:Appendix}, as they go beyond the scope of this paper.

\begin{lemma}
\label{L:diagram}

    Let $N, H, G$ groups such that $H\leq G$ and $N$ is a normal subgroup of $G$. 
    Then the inclusion homomorphism $\iota\colon H\hookrightarrow G$ induces an injective homomorphism 
\[
\kappa\colon \faktor{H}{N\cap H} \to \faktor{G}{N}.
\]
\end{lemma}

\begin{lemma}
\label{L:8groups}

        Consider the following commutative diagrams of (vertical and horizontal) short exact sequences of groups in which every square is commutative

\begin{minipage}{.48\textwidth}
\begin{equation*}
\label{E:9groups1}
\xymatrix{
  & 1 \ar[d] & 1 \ar[d] &   & \\
1 \ar[r] & U \ar[r]^-{\beta_1} \ar[d]^-{\alpha_1} & V \ar[r]^-{\pi_1} \ar[d]^-{\alpha_2} & W \ar[r] \ar[d]^-{\zeta} & 1\\
1 \ar[r] & X \ar[r]^-{\beta_2} \ar[d]^-{\mu_1} & Y \ar[r]^-{\pi_2} \ar[d]^-{\mu_2} & Z \ar[r]  & 1\\
 & R \ar[r]^-{\mu} \ar[d] & S  \ar[d] &  &  \\
  & 1  & 1  &   &
  }
\end{equation*}
\end{minipage}%
\begin{minipage}{.48\textwidth}
\begin{equation*}
\label{E:9groups2}
\xymatrix{
  &  & 1 \ar[d] & 1 \ar[d] & \\
1 \ar[r] & A \ar[r]^-{\iota_1} \ar[d]^-{\eta} & B \ar[r]^-{\rho_1} \ar[d]^-{\psi_1} & C \ar[r] \ar[d]^-{\psi_2} & 1\\
1 \ar[r] & D \ar[r]^-{\iota_2}  & E \ar[r]^-{\rho_2} \ar[d]^-{\phi_1} & F \ar[r] \ar[d]^-{\phi_2} & 1\\
 &  & G \ar[r]^-{\phi} \ar[d] & H  \ar[d] &  \\
  &   & 1  & 1  &
  }
\end{equation*}
\end{minipage}

\begin{enumerate}
    \item \begin{enumerate}
        \item If $\beta_i$ is an inclusion, for $i=1,2$, and  $\zeta$ is an isomorphism, then $\mu$ is an isomorphism.

        \item     If $\alpha_i$ is an inclusion, for $i=1,2$, and $\mu$ is an isomorphism, then $\zeta$ is an isomorphism.
    \end{enumerate}
    \label{L:8groups1}
    
    \item Suppose that, for $i=1,2$, the homomorphisms $\iota_i$ and $\psi_i$ are inclusions. 
Then $\eta$ is an isomorphism if and only if $\phi$ is.
    \label{L:8groups2}
\end{enumerate}
\end{lemma}

\proof[Proof of Theorem~\ref{T:first_main}]
From \cite[Theorem~3.1 and its proof]{ABGH} we have the following commutative diagram
\begin{equation}
\label{E:firstdiag}
\xymatrix{
1 \ar[r] & B_n[2m] \ar[r] \ar[d]^-{\psi} & B_n[m] \ar[r]^-{\tau_m} \ar[d]^-{\iota} & S_n \ar[r] \ar@{=}[d] & 1\\
1 \ar[r] & P_n \ar[r] & B_n \ar[r]^-{\tau} & S_n \ar[r] & 1}
\end{equation}
where $\tau$ is the natural surjective homomorphism which sends each braid generator $\sigma_i$ to the transposition $(i,\, i+1)$, $\tau_m$ is the restriction of $\tau$ to the subgroup $B_n[m]$, $\iota$ is the natural inclusion from the definition of congruence subgroups and $\psi$ is the restriction of $\iota$ to the subgroup $B_n[2m]$. 

Now, we consider the following diagram induced from the commutative square on the left, where the vertical arrows on this square are inclusion homomorphisms and $\psi|$ is the restriction of $\psi$, 
\begin{equation*}
\label{E:diagram}
\xymatrix{
1 \ar[r] & [P_n, P_n]\cap B_n[2m] \ar[r] \ar[d]^-{\psi|} & B_n[2m] \ar[r] \ar[d]^-{\psi} & \faktor{B_n[2m]}{[P_n, P_n]\cap B_n[2m]} \ar[r] \ar[d]^-{\overline{\psi}} & 1\\
1 \ar[r] & [P_n, P_n] \ar[r] & P_n \ar[r] & \faktor{P_n}{[P_n, P_n]} \ar[r] & 1}
\end{equation*}
From Lemma~\ref{L:diagram} the third arrow $\overline{\psi}$ on the right is also injective.
Since $\faktor{P_n}{[P_n, P_n]}$ is a free abelian group of rank $n(n-1)/2$, then $\faktor{B_n[2m]}{[P_n, P_n]\cap B_n[2m]}$ is a free abelian group of finite rank, at most $n(n-1)/2$. 
From \cite[Corollary~2.4]{ABGH} the element $A_{i,j}^m$ belongs to $B_n[2m]$, for all $1\leq i<j\leq n$, where $\{ A_{i,j} \mid 1\leq i<j\leq n \}$ is the set of Artin generators of $P_n$. 
Since $\faktor{P_n}{[P_n, P_n]}$ is generated by the set of cosets $\{ \overline{A_{i,j}} \mid 1\leq i<j\leq n \}$, it follows that $\{ \overline{A_{i,j}^m} \mid 1\leq i<j\leq n \}$ is a basis of $\faktor{B_n[2m]}{[P_n, P_n]\cap B_n[2m]}$, so it has rank $n(n-1)/2$. 
Furthermore, from the above we get 

\begin{equation}
\label{E:quot}
\frac{\left(\faktor{P_n}{[P_n, P_n]}\right)}{\left(\faktor{\overline{\psi}(B_n[2m])}{[P_n, P_n]\cap B_n[2m]}\right)} \cong (\Zz/\mm \Zz)^{n(n-1)/2}.
\end{equation}

Considering \cite[Proposition~3.1]{bloomquist2023quotients}, Arnol'd's result $B_n[2]=P_n$, and with some set theoretical equivalences, we can see that
\begin{align*}
[\PB\nn, \PB\nn] \cap \BC\nn m &= \big([\PB\nn, \PB\nn] \cap \PB\nn \big) \cap \BC\nn m \\
&= [\PB\nn, \PB\nn] \cap \big( \PB\nn \cap \BC\nn m \big) \\
&= [\PB\nn, \PB\nn] \cap \BC\nn {2m}. \\
\end{align*}

The following diagram is induced from the commutative square on the left, where the vertical arrows on this square are inclusion homomorphisms, 
\begin{equation*}
\xymatrix{
1 \ar[r] & [P_n, P_n]\cap B_n[m] \ar[r] \ar[d] & B_n[m] \ar[r] \ar[d]^-{\iota} & \faktor{B_n[m]}{[P_n, P_n]\cap B_n[m]} \ar[r] \ar[d]^-{\overline{\iota}} & 1\\
1 \ar[r] & [P_n, P_n] \ar[r] & B_n \ar[r] & \faktor{B_n}{[P_n, P_n]} \ar[r] & 1}
\end{equation*}
From Lemma~\ref{L:diagram}  the third arrow $\overline{\iota}$ on the right is also injective.
With this information, and \eqref{E:quot}, we construct the following commutative diagram 
$$
\xymatrix{
  & 1 \ar[d] & 1 \ar[d] &   & \\
1 \ar[r] & \faktor{B_n[2m]}{[P_n, P_n]\cap B_n[m]} \ar@{^{(}->}[r] \ar[d]^-{\overline{\psi}} & \faktor{B_n[m]}{[P_n, P_n]\cap B_n[m]} \ar[r] \ar[d]^-{\overline{\iota}} & S_n \ar[r] \ar@{=}[d] & 1\\
1 \ar[r] & \faktor{P_n}{[P_n, P_n]} \ar@{^{(}->}[r] \ar[d] & \faktor{B_n}{[P_n, P_n]} \ar[r] \ar[d] & S_n \ar[r]  & 1\\
 & (\Zz/\mm \Zz)^{n(n-1)/2} \ar[r]^-{\mu} \ar[d] & S  \ar[d] &  &  \\
  & 1  & 1  &   &
  }
$$
From Lemma~\ref{L:8groups} item \eqref{L:8groups1} the homomorphism $\mu$ is an isomorphism and we get the result.
\endproof

Let $n\geq 3$. Recall from \cite[Lemma~2.3]{ABGH} that the element $\sigma_i^m$ belongs to $B_n[m]$, for all $1\leq i\leq n-1$, where $\{ \sigma_i \mid 1\leq i\leq n-1 \}$ is the set of Artin generators of $B_n$.  
Although the set map $\xi\colon B_n\to B_n[m]$ defined by $\xi(\sigma_i)= \sigma_i^m$, for all $1\leq i\leq n-1$, is not a homomorphism, when $m$ is odd it induces an isomorphism on the quotient groups $\overline{\xi}\colon \faktor{\BB\nn}{[\PB\nn, \PB\nn]}\to \faktor{\BC\nn{m}}{[\PB\nn, \PB\nn] \cap \BC\nn{m}}$, as we show in the next result.

\begin{theorem}
\label{T:second_main}

Let $m$ be a positive integer and let $\nn \geq 3$. 
Consider the map 
\[
\overline{\xi}\colon \faktor{\BB\nn}{[\PB\nn, \PB\nn]}\to \faktor{\BC\nn{m}}{[\PB\nn, \PB\nn] \cap \BC\nn{m}}
\]
defined by $\overline{\xi}(\sigma_i)= \sigma_i^m$ for all $1\leq i\leq n-1$. 
If $m$ is odd, then $\overline{\xi}$ is an isomorphism. 
As a consequence, for $n\geq 3$ and $m$ odd,  $\faktor{\BC\nn{m}}{[\PB\nn, \PB\nn] \cap \BC\nn{m}}$ is a crystallographic group of dimension $n(n-1)/2$ and holonomy group $S_n$.
\end{theorem}

\proof
Suppose that $n\geq 3$ and $m$ is an odd positive integer and consider the map 
\[
\overline{\xi}\colon \faktor{\BB\nn}{[\PB\nn, \PB\nn]}\to \faktor{\BC\nn{m}}{[\PB\nn, \PB\nn] \cap \BC\nn{m}}
\]
defined by $\overline{\xi}(\sigma_i)= \sigma_i^m$ for all $1\leq i\leq n-1$.
To show that $\overline{\xi}$ is a homomorphism it is enough to verify that Artin's relations are preserved by $\overline{\xi}$. 

Let $1\leq i,j\leq n$ such that $|i-j|\geq 2$. 
From Artin's relation $\sigma_i\sigma_j=\sigma_j\sigma_i$ we obtain $\sigma_i^m\sigma_j^m=\sigma_j^m\sigma_i^m$ in $B_n[m]$, which is then preserved  by~$\overline{\xi}$. 

Let $1\leq i\leq n-2$. The equality $\sigma_i^m\sigma_{i+1}^m\sigma_i^m\sigma_{i+1}^{-m}\sigma_i^{-m}\sigma_{i+1}^{-m}=1$ is valid in $\faktor{\BC\nn{m}}{[\PB\nn, \PB\nn] \cap \BC\nn{m}}$. 
In fact, suppose that $m=2k+1$, then from the action of conjugation in $\faktor{B_n}{[P_n, P_n]}$ described in Subsection~\ref{Subsec:braids} we have 
\begin{align*}
\sigma_i^m\sigma_{i+1}^m\sigma_i^m\sigma_{i+1}^{-m}\sigma_i^{-m}\sigma_{i+1}^{-m} &= 
A_{i,i+1}^k\sigma_i A_{i+1,i+2}^k\sigma_{i+1}A_{i,i+1}^k\sigma_i \sigma_{i+1}^{-1}A_{i+1,i+2}^{-k} \sigma_i^{-1}A_{i,i+1}^{-k} \sigma_{i+1}^{-1}A_{i+1,i+2}^{-k}\\
&= A_{i,i+1}^k A_{i,i+2}^k A_{i+1,i+2}^k A_{i,i+1}^{-k} A_{i,i+2}^{-k} A_{i+1,i+2}^{-k} \\
&= 1 \in \faktor{B_n}{[P_n, P_n]}. 
\end{align*}
From Theorem~\ref{T:first_main} the homomorphism 
\[
\overline{\iota}\colon \faktor{B_n[m]}{[P_n, P_n]\cap B_n[m]}\to \faktor{B_n}{[P_n, P_n]}
\]
is injective, then 
$\sigma_i^m\sigma_{i+1}^m\sigma_i^m\sigma_{i+1}^{-m}\sigma_i^{-m}\sigma_{i+1}^{-m}=1$ in $\faktor{\BC\nn{m}}{[\PB\nn, \PB\nn] \cap \BC\nn{m}}$.

Now, consider the following commutative diagram of short exact sequences
\[
\xymatrix{
1 \ar[r] & \faktor{P_n}{[P_n, P_n]} \ar[r] \ar[d]^-{\overline{\xi}|} & \faktor{B_n}{[P_n, P_n]} \ar[r] \ar[d]^-{\overline{\xi}} & S_n \ar[r] \ar@{=}[d] & 1\\
1 \ar[r] & \faktor{B_n[2m]}{[P_n, P_n]}\cap B_n[m] \ar[r]  & \faktor{B_n[m]}{[P_n, P_n]\cap B_n[m]} \ar[r]  & S_n \ar[r]  & 1  }
\]
Recall from the proof of Theorem~\ref{T:first_main} that the free abelian groups $\faktor{B_n[2m]}{[P_n, P_n]\cap B_n[2m]}$ and  $\faktor{P_n}{[P_n, P_n]}$ of rank $n(n-1)/2$ have a basis $\{ \overline{A_{i,j}^m} \mid 1\leq i<j\leq n \}$ and $\{ \overline{A_{i,j}} \mid 1\leq i<j\leq n \}$, respectively.
Since 
\[
\overline{\xi}|\colon \faktor{P_n}{[P_n, P_n]}\to \faktor{B_n[2m]}{[P_n, P_n]\cap B_n[2m]}
\]
is a homomorphism such that $\overline{\xi}|(A_{i,j})= A_{i,j}^m$, for all $1\leq i<j\leq n$, then it is an isomorphism. 
Therefore, from the five lemma $\overline{\xi}$ is an isomorphism. 

The last part follows from the result on the crystallographic braid group $\faktor{B_n}{[P_n, P_n]}$, see \cite[Proposition~1]{Goncalves-Guaschi-Ocampo:2017}.
\endproof

A group $G$ is called \textit{co-Hopfian} if it is not isomorphic to any of its proper subgroups, or equivalently if  every injective homomorphism $\phi\colon G\to G$ is surjective. 
It is known that the braid group $B_n$ is not co-Hopfian. However, for $n\geq 4$, the quotient by its center is co-Hopfian, see Bell-Margalit \cite{Bell-Margalit:2006}. 

\begin{corollary}
    Let $n\geq 3$. The crystallographic braid group $\faktor{B_n}{[P_n, P_n]}$ is not co-Hopfian.
\end{corollary}

\proof
It follows from Theorems~\ref{T:first_main} and~\ref{T:second_main}. 
\endproof

\appendix

\section{Technical proofs}\label{S:Appendix}

In this section we give proofs for Lemmas~\ref{L:diagram} and~\ref{L:8groups}.

\proof[Proof of Lemma \ref{L:diagram}]
From the hypothesis we have the following commutative square of inclusions
$$
\xymatrix{
N\cap H \ar[r] \ar[d]^-{\iota|} & H \ar[d]^-{\iota}\\
N \ar[r] & G }
$$
that induces the following commutative diagram of short exact sequences and the existence of the homomorphism $\kappa\colon \faktor{H}{N\cap H} \to \faktor{G}{N}$
\begin{equation}
\label{E:diaglemma}
\xymatrix{
1 \ar[r] & N\cap H \ar[r] \ar[d]^-{\iota|} & H \ar[r]^-{\pi} \ar[d]^-{\iota} & \faktor{H}{N\cap H} \ar[r] \ar[d]^-{\kappa} & 1\\
1 \ar[r] & N \ar[r] & G \ar[r]^-{\rho} & \faktor{G}{N} \ar[r] & 1}
\end{equation}
Now we use diagram chasing to prove that $\kappa$ is injective. 
Let $a\in \faktor{H}{N\cap H}$ such that $\kappa(a)=1$. 
Since $\pi$ is surjective, then there exists $h\in H$ such that $\pi(h)=a$. 
From the commutativity of the diagram we get $\rho(\iota(h))=1$, i.e. $\rho(h)=1$. 
So, $h$ belongs to the kernel of $\rho$ that is equal to $N$. 
Hence, $h\in N\cap H$ and from the short exact sequence on the top of \eqref{E:diaglemma} we conclude $a=1$. 
\endproof

\proof[Proof of Lemma \ref{L:8groups}]
We use diagram chasing, the commutativity of the squares and the definition of short exact sequences to complete some steps of the proof. 

\begin{enumerate}
    \item 
    \begin{enumerate}
        \item   Suppose that $\beta_i$ is an inclusion, for $i=1,2$. Then
    from Lemma~\ref{L:diagram}  the homomorphisms $\mu$ is injective.
Suppose that $\zeta$ is an isomorphism. We will prove that $\mu$ is surjective. 
    Let $s\in S$, then there is $y\in Y$ such that $\mu_2(y)=s$. Let $z=\pi_2(y)$. 
    Since $\zeta$ is an isomorphism, there exists $w\in W$ such that $\zeta(w)=z$ and from the surjectivity of $\pi_1$ there is $v\in V$ such that $\pi_1(v)=w$. 
    Then from the commutativity of the diagram we have $\pi_2(v^{-1}y)=1$ that means $v^{-1}y\in X$. 
    Again, from the commutativity of the diagram, we conclude that the element $r=\mu_1(v^{-1}y)$ is such that $\mu(r)=s$. 
    Then $\mu$ is surjective. 
    
    \item The proof is similar to the first item.
    \end{enumerate}

\item From the five lemma we get that $\eta$ is injective and that $\phi$ is surjective.

Suppose that $\eta$ is an isomorphism. 
We  prove that $\phi$ is injective. 
Let $g\in G$ such that $\phi(g)=1$. 
Then there is an $e\in E$ such that $\phi_1(e)=g$. Let $f=\rho_2(e)$.
By using the commutativity of the squares we have $\phi_2( f )=1$, so we have $f\in C$. 
From the surjection $\rho_1$, there is a $b\in B$ such that $\rho_1(b)=f$. 
Since $\psi_1$ is an inclusion and $\rho_2(e)=f$, we have $b^{-1}e\in Ker(\rho_2)=D$. 
We know that  $\eta$ is an inclusion, and from the hypothesis we know it is an isomorphism. We can conclude that it is the identity map. 
Hence, $b^{-1}e\in A$ and $e$ belongs to $B$, so $g=1$. 
We conclude that $\phi$ is injective.

Similarly we prove that $\eta$ is surjective if $\phi$ is an isomorphism.
\end{enumerate}
\endproof

\begin{remark}
 We note that in this paper we do not use Lemma~\ref{L:8groups} item \eqref{L:8groups2}. However, it will be useful in~\cite{BDMOS:2023}.
\end{remark}

\bibliography{braid_groups.bib}{}
\bibliographystyle{alpha}

\end{document}